\documentclass[12pt,oneside]{article}
\usepackage[cp866]{inputenc}   
\usepackage{color}
\usepackage{amsmath}
\usepackage{amssymb}
\usepackage{amsthm}
\usepackage[english]{babel}

\newtheorem{theo}{\bf Theorem}
\newtheorem{lemma}{\bf Lemma}
\newcommand {\bproof }{{\par\medskip\noindent \bf Proof. }}

\def\wbull{\hfill\vrule height .9ex width .8ex depth -.1ex}

\title{{\bf {\bf \normalsize Alexander L. Gavrilyuk\footnote{
~e-mail:~\texttt{alexander.gavriliouk@gmail.com}\newline
Institute of Mathematics and Mechanics of Russian Academy of Sciences, Ekaterinburg, Russia.},~
Ivan Y. Mogilnykh\footnote{
~e-mail:~\texttt{ivmog84@gmail.com}\newline
Institute of Mathematics of Russian Academy of Sciences, Novosibirsk, Russia.}\\}
\normalsize \large CAMERON -- LIEBLER LINE CLASSES IN ${\rm PG}(n,4)$, $n\geqslant 3$}}
\begin{document}

\maketitle
\section{Introduction}
\hspace{5mm}

A {\it Cameron -- Liebler line class} ${\cal L}$ with parameter $x$ is
a set of lines of projective geometry $PG(3,q)$ such that each
line of ${\cal L}$ meets exactly $x(q+1)+q^2-1$ lines of ${\cal
L}$ and each line that is not from ${\cal L}$ meets exactly
$x(q+1)$ lines of ${\cal L}$ (there are several equivalent
definitions of Cameron -- Liebler line classes, see Section 2). These
classes appeared in connection with an attempt by Cameron and
Liebler  \cite{CameronLiebler} to classify collineation groups
of $PG(n,q)$, $n\geqslant 3$, that have equally many orbits on lines and on points.

The following line classes (and their complementary sets)
are examples of Cameron -- Liebler line classes: 
\smallskip

-- the set of all lines in $PG(3,q)$,
\smallskip

-- the set of all lines in a given plane of $PG(3,q)$,
\smallskip

-- the set of all lines through a point,
\smallskip

-- for a non-incident point -- hyperplane pair $(P,\pi)$, the set of all lines through $P$ or in $\pi$.
\smallskip

Cameron and Liebler conjectured \cite{CameronLiebler} that, apart
from these examples, there are no Cameron -- Liebler line classes.
The counterexamples were constructed by Drudge \cite{Drudge}
(in $PG(3,3)$ with $x=5$), by Bruen and Drudge \cite{BruenDrudge}
(for odd $q$, in $PG(3,q)$ with $x=(q^2+1)/2$),
by Govaerts and Penttila \cite{GovaertsPenttila} (in $PG(3,4)$ with $x=7$),
and recently by Rodgers \cite{Rodgers} (for some odd $q$, in $PG(3,q)$ with $x=(q^2-1)/2$).
A complete classification of Cameron -- Liebler line classes in $PG(3,3)$
was obtained by Drudge \cite{Drudge}. Cameron -- Liebler line classes in $PG(3,4)$
were studied by Govaerts and Penttila \cite{GovaertsPenttila}.
However, they left two open cases with parameter $x\in \{6,8\}$.

In this paper, we show the non-existence of Cameron -- Liebler line
classes with parameter $x\in \{6,8\}$ in $PG(3,4)$ and give a new
proof of non-existence of those with parameter $x\in \{4,5\}$ in $PG(3,4)$
previously established in \cite{GovaertsPenttila}. Further,
we prove the uniqueness of Cameron -- Liebler line class with $x=7$
in $PG(3,4)$ discovered in \cite{GovaertsPenttila}. Finally,
following the approach by Drudge \cite{DrudgeThesis}, we
obtain a complete classification of Cameron -- Liebler line classes
in $PG(n,4)$, $n\geqslant 3$ (for the precise definition of those in $PG(n,q)$,
see Section 5).

In our proof we consider a Cameron -- Liebler line class as a
subset of the vertex set of the Grassmann graph. Recall that the
{\it Grassmann graph} $G_q(n,e)$ is a graph whose vertex set
consists of all $e$-dimensional subspaces of an $n$-dimensional
vector space over a finite field of order $q$; two $e$-subspaces
are adjacent if and only if their intersection has dimension
$e-1$. In case of $e=2$ the graph $G_q(n+1,2)$ can be viewed as a
graph on lines in $PG(n,q)$ with two distinct lines adjacent if
and only if they meet. There is a natural correspondence between
Cameron -- Liebler line classes in $PG(3,q)$ and completely
regular subsets of the vertex set of $G_q(4,2)$ with strength 0
and covering radius 1 \cite{Vanhove}.


The paper is organized as follows. In Section 2, we recall some definitions
and certain properties of Cameron -- Liebler line classes in $PG(n,q)$,
rewrite them in terms of the Grassmann graphs and formulate a 
new necessary condition. In Section 3, we describe properties of
putative Cameron -- Liebler line classes in $PG(3,4)$ and then we
obtain a contradiction by simple counting arguments. Section 4 is
devoted to the uniqueness of Cameron -- Liebler line class with
$x=7$ in $PG(3,4)$. In Section 5 we obtain a classification of
Cameron -- Liebler line classes in $PG(n,4)$, $n\geqslant 3$ as a
consequence of results from the previous sections.


\section{Properties of Cameron -- Liebler line classes}

Following Penttila \cite{Penttila}, for a point $P$ of $PG(3,q)$,
we denote the set of all lines through $P$ by ${\rm Star}(P)$, and, for
a hyperplane $\pi$ of $PG(3,q)$, the set of all lines in $\pi$ by
${\rm line}(\pi)$. Both types of line sets will be also reffered to as a
{\it clique} in $PG(3,q)$. For a hyperplane $\pi$ and a point
$P\in \pi$, a {\it pencil} ${\rm pen}(P,\pi)$ is the set of all
lines in $\pi$ through $P$.

For a set of lines ${\cal L}$ in $PG(3,q)$, let ${\overline {\cal L}}$
denote a complementary set of lines, and $\chi_{\cal L}$ denote
the characteristic function of ${\cal L}$.

\begin{lemma}\label{Definitions} $(\cite{Penttila})$
The following conditions are equivalent.

\quad $(1)$ There exists an integer $x$ such that, for every line
$l$,
$$|\{m\in {\cal L}\setminus \{l\}:~m~meets~l\}| = (q+1)x + (q^2-1)\chi_{\cal L}(l).$$

\quad $(2)$ There exists an integer $x$ such that, for every
incident point-plane pair $(P,\pi)$,
$$|{\rm Star}(P)\cap {\cal L}| + |{\rm line}(\pi)\cap {\cal L}| = x + (q+1)|{\rm pen}(P,\pi)\cap {\cal L}|.$$

\quad $(3)$ There exists an integer $x$ such that, for every pair
of skew lines $l$ and $m$,
$$|\{n\in {\cal L}\setminus \{l,m\}:~n~meets~l~and~n~meets~m\}| = x + q(\chi_{\cal L}(l)+\chi_{\cal L}(m)).$$
\end{lemma}

A set of lines ${\cal L}$ is called a {\it Cameron -- Liebler
line class} in $PG(3,q)$ if one of the conditions in Lemma
\ref{Definitions} is satisfied. 
We note that the number $x$ in
each of these conditions is the same and is called the
{\it parameter} of the Cameron -- Liebler line class. In
\cite{CameronLiebler} it is shown that $|{\cal L}|=x(q^2+q+1)$
holds so that $x\in\{0,1,2,\dots,q^2+1\}$, and ${\overline {\cal
L}}$ is a Cameron -- Liebler line class with parameter $q^2+1-x$.
\smallskip

A graph $G$ defined on the set of lines of $PG(3,q)$
with two distinct lines adjacent if and only if they meet
is the well-known Grassmann graph $G_q(4,2)$.
The Grassmann graph $G_q(n,e)$ is distance-transitive
and has diameter ${\rm min}\{n-e,e\}$. A detailed discussion of
the Grassmann graphs is contained in \cite{BCN}. Here we recall
some properties necessary for our partial case (Lemma
\ref{GrassmannProperties} below, see \cite[Chapter 9.3]{BCN} for
its proof).

If $X$ is a subset of the vertex set of $G:=G_q(4,2)$
then in order to shorten the notation we write $X$ for the graph
induced by $G$ on $X$. For vertices $v,u\in G$, we define the {\it
neighborhood} $G(v):=\{w\in G|~w\sim v\}$ of $v$, the {\it second
neighborhood} $G_2(v):=\{w\in G\setminus \{v\}|~w\not\sim v\}$,
and $G(u,v):=G(u)\cap G(v)$.

For an integer $\alpha\geqslant 1$, the $\alpha$-{\it clique extension}
of a graph ${\overline H}$ is the graph $H$ obtained from
${\overline H}$ by replacing each vertex ${\overline u}\in
{\overline H}$ with a clique $U$ with $\alpha$ vertices, where,
for any ${\overline u},{\overline w}\in {\overline H}$, $u\in U$
and $w\in W$, ${\overline u}$ and ${\overline w}$ are adjacent if
and only if $u$ and $w$ are adjacent. By the $n\times m$-{\it
grid}, we mean the Cartesian product of two cliques on $n$ and $m$
vertices.

\begin{lemma}\label{GrassmannProperties}
The following holds.

\quad $(1)$ For every vertex $v\in G$, the graph $G(v)$
is the $q$-clique extension of $(q+1)\times (q+1)$-grid
$($so that $G$ is regular with valency $q(q+1)^2$$)$.

\quad $(2)$ For every pair of non-adjacent vertices $u,v\in G$,
the graph $G(u,v)$ is the $(q+1)\times (q+1)$-grid.
\end{lemma}

It is easily seen that there are exactly two different sets of
maximal cliques  in $G$, say $\Lambda_1$ and $\Lambda_2$, with
each maximal clique having the same size $q(q+1)+1$. We define
$\Lambda_1$ to be the set of all cliques that correspond to
${\rm line}(\pi)$, for all hyperplanes $\pi$ in $PG(3,q)$, and
$\Lambda_2$ to be the set of all cliques that correspond to
${\rm Star}(P)$, for all points $P$ in $PG(3,q)$. Further, for every
pair of cliques $L,L'\in \Lambda_i$, we have $|L\cap L'|=1$,
while, for every pair of cliques $L\in \Lambda_1$, $L^{*}\in
\Lambda_2$, we have $|L\cap L^{*}|\in \{0,q+1\}$, and, for every
pair of adjacent vertices $u,v\in G$, there is a unique pair of
cliques $L\in \Lambda_1$ and $L^{*}\in \Lambda_2$ such that
$u,v\in L\cap L^{*}$.
\smallskip

Now let ${\cal L}$ be a Cameron -- Liebler line class in $PG(3,q)$ with parameter $x$.
Consider ${\cal L}$ as a subset of vertices of $G_{q}(4,2)$.
Clearly, a partition of the vertex set of $G$
into two parts, ${\cal L}$ and ${\overline {\cal L}}$, is {\it equitable}
with the following quotient matrix
\[ P:=\bordermatrix{
                     & {\cal L} & {\overline {\cal L}} \cr
{\cal L}             & (q+1)x+q^2-1 & q(q+1)^2-(q+1)x-q^2+1 \cr
{\overline {\cal L}} & (q+1)x       & q(q+1)^2-(q+1)x},\] which
means that every vertex from a part $A$ has exactly $p_{A,B}$
neighbors in a part $B$.

\begin{lemma}\label{DefinitionsInGraphs}
The following conditions are equivalent.

\quad $(1)$ for every vertex $v\in G$,
\begin{equation}\label{Valency}
|G(v)\cap {\cal L}| = (q+1)x + (q^2-1)\chi_{\cal L}(v).
\end{equation}

\quad $(2)$ for every pair of cliques $L\in \Lambda_1$, $L^{*}\in \Lambda_2$
such that $|L\cap L^*|=q+1$,
\begin{equation}\label{StarLinePencil}
|L\cap {\cal L}|+|L^{*}\cap {\cal L}|=x+(q+1)|L\cap L^{*}\cap {\cal L}|.
\end{equation}

\quad $(3)$ for every pair of non-adjacent vertices $u,v\in G$,
\begin{equation}\label{Mu}
|G(u,v)\cap {\cal L}| = x + q(\chi_{\cal L}(u)+\chi_{\cal L}(v)).
\end{equation}

\quad $(4)$ for every pair of adjacent vertices $u,v\in G$,
\begin{equation}\label{B1}
|G(v)\cap G_2(u)\cap {\cal L}|=q (x + q\chi_{\cal L}(v) - |L\cap L^{*}\cap {\cal L}|),
\end{equation}
where $L\in \Lambda_1$, $L^{*}\in \Lambda_2$, and $u,v\in L\cap L^{*}$.
\end{lemma}
\bproof The first three statements are equivalent to Lemma
\ref{Definitions}. Let us show (4). By Eq. (\ref{Valency}), for a
vertex $v$, we have $|G(v)\cap {\cal L}| = (q+1)x +
(q^2-1)\chi_{\cal L}(v)$, while $G(v)$ contains exactly $|L\cap
{\cal L}|+|L^{*}\cap {\cal L}|-|L\cap L^{*}\cap {\cal
L}|-\chi_{\cal L}(v)$ vertices from ${\cal L}\cap (L\cup
L^{*}\setminus \{v\})\subset G(u)\cup \{u\}$. Therefore,
$$|G(v)\cap G_2(u)\cap {\cal L}|=
(q+1)x + (q^2-1)\chi_{\cal L}(v)-
(|L\cap {\cal L}|+|L^{*}\cap {\cal L}|-|L\cap L^{*}\cap {\cal L}|-\chi_{\cal L}(v))$$

Taking into account (\ref{StarLinePencil}), we obtain
$|G(v)\cap G_2(u)\cap {\cal L}|=
q(x + q\chi_{\cal L}(v)-|L\cap L^{*}\cap {\cal L}|).$

Now suppose (4) holds. Summing Eq. (\ref{B1}) over  all vertices
$u\in G(v)$, we obtain

$$\sum_{u\in G(v)}|G(v)\cap G_2(u)\cap {\cal L}|=
q^2(q+1)^2x + q^3(q+1)^2\chi_{\cal L}(v)-q^2\sum_{L\in \Lambda_1,L^{*}\in \Lambda_2:~v\in L\cap L^{*}}
|L\cap L^{*}\cap {\cal L}|.$$

Note that
$$\sum_{L\in \Lambda_1,L^{*}\in \Lambda_2:~v\in L\cap L^{*}}|L\cap L^{*}\cap {\cal L}|=
|G(v)\cap {\cal L}|+(q+1)^2\chi_{\cal L}(v).$$

On the other hand, for every vertex $u\in G(v)$,
$|G(v)\cap G_2(u)|=q^3$ holds, hence
$$\sum_{u\in G(v)}|G(v)\cap G_2(u)\cap {\cal L}|=q^3|G(v)\cap {\cal L}|,$$
which yields that (1) holds. This proves the lemma.\wbull
\medskip

Many of the previous results on Cameron--Liebler line classes were obtained by Drudge's approach. This approach relates the
lines of a Cameron--Liebler line class that belong to a
clique in $PG(3,q)$ to a blocking set in a projective plane
$PG(2,q)$. We recall that a {\it t-fold blocking set} in $PG(2,q)$
is a set of points that intersects every line in at least $t$
points. A 1-fold blocking set is called just a blocking set, and
it is called {\it trivial} if it contains a line. It is also
reasonable to recall that every clique in $PG(3,q)$ and its lines
correspond to a projective plane and its points respectively,
while pencils in the clique correspond to lines in the plane. In
this way, for a clique $L$, some restrictions on $|L\cap {\cal
L}|$ may be obtained from the study of blocking sets, \cite{GovaertsPenttila}.
\smallskip

Our main results in the next sections rely on the study of a
possible distribution of lines from a Cameron -- Liebler line
class ${\cal L}$ in the set $\cup_{P\in u}\cup_{\pi\ni u} {\rm
pen}(P,\pi)$ for a line $u$. In terms of the Grassmann graph $G$,
for a vertex $u$, we consider the graph $G(u)$, which is the
$q$-clique extension of $(q+1)\times (q+1)$-grid by Lemma
\ref{GrassmannProperties}(1). 
The intersection of every pair of cliques $L\in \Lambda_1$, $L^*\in \Lambda_2$
that contain $u$ is a $q$-clique in $G(u)$, and it corresponds to ${\rm pen}(P,\pi)\setminus
\{u\}$ for some point $P\in u$ and hyperplane $\pi\ni u$. Further,
define a square matrix ${\cal T}(u)$ of size $q+1$,
$${\cal T}(u):=
\left (
\begin{array}{cccc}
t_{1,1} & t_{1,2} & \ldots &  t_{1,q+1}\\
\ldots  & \ldots  & \ldots &  \ldots\\
t_{q+1,1} & t_{q+1,2} & \ldots &  t_{q+1,q+1}
\end{array}
\right ),$$
whose elements are equal to $|{\cal L}\cap {\rm pen(P,\pi)}\setminus \{u\}|$,
so that, without the loss of generality, the rows (resp. columns) of ${\cal T}(u)$ correspond
to points on $u$ (resp. hyperplanes containing $u$). Clearly, we do not distinguish
between matrices obtained by permutations of rows or columns as well as transposition.

For a line $u$, we call the matrix ${\cal T}(u)$ the {\it pattern} with respect to $u$
(the {\it pattern} for short). Then a line $u$ {\it has} pattern ${\cal T}(u)$.

\begin{lemma}\label{NewCondition}
For every vertex $u\in G$, the pattern ${\cal T}(u)$ satisfies the following properties:

\quad $(1)$ $0\leqslant t_{ij}\leqslant q$ for all $i,j\in \{1,\dots,q+1\}$$;$

\quad $(2)$ ${\displaystyle \sum_{i,j=1}^{q+1}t_{ij}=x(q+1)+\chi_{\cal L}(u) (q^2-1)}$$;$

\quad $(3)$ ${\displaystyle \sum_{j=1}^{q+1}t_{kj}+\sum_{i=1}^{q+1}t_{il}=x+(q+1)(t_{kl}+\chi_{\cal L}(u))}$,
for all $k,l\in \{1,\dots,q+1\}$$;$

\quad $(4)$ ${\displaystyle \sum_{i,j=1}^{q+1}t_{ij}^2=(x-\chi_{\cal L}(u))^2+q(x-\chi_{\cal L}(u))+\chi_{\cal L}(u) q^2(q+1)}$.
\end{lemma}
\bproof
Let $u$ be a vertex of $G$.
Let $L_1,L_2,\dots,L_{q+1}$ ($L^1,L^2,\dots,L^{q+1}$, respectively)
be the set of all maximal cliques from $\Lambda_1$ (from $\Lambda_2$,
respectively) containing $u$. Then, by definition, $t_{ij}:=|L_i\cap L^j\cap {\cal L}\setminus \{u\}|$
(so that $|L_i\cap {\cal L}|=\chi_{\cal L}(u)+\sum_{j} t_{ij}$ and
$|L^j\cap {\cal L}|=\chi_{\cal L}(u)+\sum_{i} t_{ij}$).
Now Statement (1) is clear, and Statements (2) and (3) follow from
Statements (1) and (2) of Lemma \ref{DefinitionsInGraphs}, respectively.

Let us prove (4). Recall that $|{\cal L}|=x(q^2+q+1)$. By Lemma
\ref{DefinitionsInGraphs}(\ref{Mu}), this yields that there are
exactly
$$(q(\chi_{\cal L}(u)+1)+x)|G_2(u)\cap {\cal L}|$$
edges $\{v,w\}$ such that $v\in G(u)\cap {\cal L}$, $w\in
G_2(u)\cap {\cal L}$, where
$$|G_2(u)\cap {\cal L}|=x(q^2+q+1)-x(q+1)-\chi_{\cal L}(u)(q^2-1)-\chi=q^2(x-\chi_{\cal L}(u)).$$

On the other hand, by Lemma \ref{DefinitionsInGraphs}(\ref{B1}),
the same number must be equal to
$$\sum_{v\in G(u)\cap {\cal L}}|G(v)\cap G_2(u)\cap {\cal L}|=
\sum_{i,j=1}^{q+1}t_{ij}q(x+q-t_{ij}-\chi_{\cal L}(u)),$$ which
proves the lemma after some straightforward calculations.\wbull
\medskip

Statement (4) of Lemma \ref{NewCondition} is a new existence condition
for Cameron -- Liebler line classes ${\cal L}$ in $PG(3,q)$, which follow,
in fact, from combining Statements (3) and (4) of Lemma \ref{DefinitionsInGraphs}.

In general, a pattern w.r.t. a line is determined by any pair of
its row and column. Hence, all the possible patterns w.r.t. a line
may be obtained by enumerating all the row-column pairs and
checking the properties (1)--(4) of Lemma \ref{NewCondition}.
Sometimes it is enough to show the non-existence of a putative
Cameron -- Liebler line class in the sense that the set of
possible patterns turns out to be empty, see Section 3.

The following remark may be 
useful in the study of structure of a Cameron -- Liebler line class 
(although we do not involve it in Section 4, we think that it could be
used in an alternative proof of the uniqueness of Cameron --
Liebler line class with $x=7$ in $PG(3,4)$).

Let $u,v$ be a pair of vertices of $G$ at distance 2.
Consider the subgraph $G(u,v)$, and recall that, by Lemma \ref{GrassmannProperties},
it is the $(q+1)\times (q+1)$-grid.
Let $L_1,L_2,\dots,L_{q+1}$ ($L^1,L^2,\dots,L^{q+1}$, respectively)
be the set of all maximal cliques from $\Lambda_1$ (from $\Lambda_2$,
respectively) containing $u$, and $K_1,K_2,\dots,K_{q+1}$ ($K^1,K^2,\dots,K^{q+1}$,
respectively) be the set of all maximal cliques from $\Lambda_1$ (from $\Lambda_2$,
respectively) containing $v$.

Let the cliques $L_1,L_2,\ldots,L_{q+1}$, $L^1,L^2,\ldots,L^{q+1}$,
$K_1,K_2,\ldots,K_{q+1}$, $K^1,K^2,\ldots,K^{q+1}$ be ordered
so that each of the maximal $(q+1)$-cliques of $G(u,v)$
is the intersection of the type $L_i\cap K^i$ or $L^j\cap K_j$,
$i,j\in\{1,\ldots,q+1\}$.

Set $m_i:=|L_i\cap K^i\cap {\cal L}|$, $n_j:=|L^j\cap K_j\cap {\cal L}|$,
$i,j\in\{1,\ldots,q+1\}$.
Then, by Lemma \ref{DefinitionsInGraphs}, $m_i,n_j$ are
non-negative integers, each of them is at most $q+1$, and the following holds:
$$|L_i\cap {\cal L}|+|K^i\cap {\cal L}|=x+(q+1)m_i,~
|L^j\cap {\cal L}|+|K_j\cap {\cal L}|=x+(q+1)n_j,$$
$$\sum m_i=\sum n_j=x+q(\chi_{\cal L}(u)+\chi_{\cal L}(v))=|G(u,v)\cap {\cal L}|.$$

We note that the set $G(u,v)$ can be naturally associated with a
$(q+1)\times (q+1)$-matrix $M$, whose entry equals 1 or 0 whenever
the corresponding vertex of $G(u,v)$ belongs to ${\cal L}$ or
${\overline {\cal L}}$. We may assume that the $i$th row of $M$
contains exactly $m_i$ ones, while the $j$th column of $M$
contains exactly $n_j$ ones.

The question of existence of such a matrix $M$ is quite
well-studied and leads us to the so-called {\it Ryser} classes of
$(0,1)$-matrices with given row and column sums, the non-emptiness
of which can be easily settled, see \cite[Theorem 1.2.8]{Sachkov}.

%
%



\section{Non-existence of Cameron -- Liebler line classes with $x\in\{4,5,6,8\}$ in $PG(3,4)$}

In this section, we show that there are no Cameron -- Liebler line
classes in $PG(3,4)$ with parameter $x\in \{6,8\}$ or,
equivalently, the Grassmann graph $G_4(4,2)$ does not admit
equitable partitions with quotient matrices
$$\left ( \begin{array}{cc}
45 & 55 \\
30 & 70
\end{array}\right ), {\rm~and~}\left ( \begin{array}{cc}
55 & 45 \\
40 & 60
\end{array}\right ).$$

The non-existence of Cameron -- Liebler line classes with parameter
$x\in \{4,5\}$ was shown in \cite[Theorem 1.3]{GovaertsPenttila}.
Our arguments can also cover these cases.

\begin{theo}\label{NonExistence}
There are no Cameron -- Liebler line classes with parameter $x\in
\{4,5,6,8\}$ in $PG(3,4)$.
\end{theo}
\bproof Let $G:=G_4(4,2)$ and ${\cal L}$ be a Cameron -- Liebler line
class with parameter 6 in $PG(3,4)$.
First of all, we need the following lemma \cite[Theorem 1.3(3)]{GovaertsPenttila}.

\begin{lemma}\label{Penttila6}
If ${\cal L}$ is a Cameron -- Liebler line class with parameter $6$ in $PG(3,4)$,
then every maximal clique of $G$ intersects ${\cal L}$ in $3~{\rm mod}~5$ vertices.
Moreover, for each $\alpha\in \{3,8,13,18\}$, there exists a maximal clique containing
exactly $\alpha$ vertices of ${\cal L}$.
\end{lemma}

We follow the notation used in the proof of Lemma
\ref{NewCondition}. Let $u$ be a vertex from ${\overline {\cal
L}}$.
By Lemma \ref{Penttila6}, without loss of generality,
we may assume that $|L_1\cap {\cal L}|=\sum_{j=1}^{5}t_{1j}=3$,
and consider the following three possibilities for $t_{1j}$, $j=1,\dots,5$:
\begin{itemize}
\item $t_{11}=3$, $t_{12}=\dots=t_{15}=0$,
\item $t_{11}=2$, $t_{12}=1$, $t_{13}=t_{14}=t_{15}=0$,
\item $t_{11}=t_{12}=t_{13}=1$, $t_{14}=t_{15}=0$,
\end{itemize}
and the following three possibilities for $|L_k\cap {\cal L}|=\sum_{j=1}^{5}t_{kj}$, $k=1,\dots,5$:
\begin{itemize}
\item $|L_k\cap {\cal L}|=3$, $k=1,\dots,4$, and $|L_5\cap {\cal L}|=18$,
\item $|L_k\cap {\cal L}|=3$, $k=1,2,3$, and $|L_4\cap {\cal L}|=8$, $|L_5\cap {\cal L}|=13$,
\item $|L_k\cap {\cal L}|=3$, $k=1,2$, and $|L_k\cap {\cal L}|=8$, $k=3,4,5$.
\end{itemize}

Given the values of $t_{1j}$, $j=1,\dots,5$, and $|L_k\cap {\cal
L}|$, $k=1,\dots,5$, Lemma \ref{NewCondition}(3) allows us to
determine the remaining elements of pattern ${\cal
T}(u):=(t_{ij})_{5\times 5}$. Indeed, the numbers $|L^l\cap {\cal
L}|=\sum_{i=1}^{5}t_{il}$, $l=1,\dots,5$, are calculated from
Lemma \ref{NewCondition}(3) when $k=1$:
$$\sum_{j=1}^{5}t_{1j}+\sum_{i=1}^{5}t_{il}=6+5t_{1l},$$
and, further, for every pair of indices $k,l\in \{1,\dots,5\}$,
the value of $t_{kl}$ is determined by $\sum_{j=1}^{5}t_{kj}$
and $\sum_{i=1}^{5}t_{il}$.

Taking into account Statements (1), (2), and (3) of Lemma
\ref{NewCondition}, we obtain only 6 admissible variants for
${\cal T}(u)$. However, two pairs of them are equivalent under the
action of automorphism of $G$ that interchanges $\Lambda_1$ and
$\Lambda_2$ (see \cite[Theorem 9.3.1]{BCN}). Therefore, we have
the following candidates for ${\cal T}(u)$:
$$ T_1 = \left(%
\begin{array}{ccccc}
  1 & 1 & 1 & 0 & 0 \\
  1 & 1 & 1 & 0 & 0 \\
  1 & 1 & 1 & 0 & 0 \\
  1 & 1 & 1 & 0 & 0 \\
  4 & 4 & 4 & 3 & 3 \\
\end{array}%
\right),
T_2 = \left(%
\begin{array}{ccccc}
  2 & 1 & 0 & 0 & 0 \\
  2 & 1 & 0 & 0 & 0 \\
  2 & 1 & 0 & 0 & 0 \\
  3 & 2 & 1 & 1 & 1 \\
  4 & 3 & 2 & 2 & 2 \\
\end{array}%
\right),
T_3 = \left(%
\begin{array}{ccccc}
  1 & 1 & 1 & 0 & 0 \\
  1 & 1 & 1 & 0 & 0 \\
  1 & 1 & 1 & 0 & 0 \\
  2 & 2 & 2 & 1 & 1 \\
  3 & 3 & 3 & 2 & 2 \\
\end{array}%
\right),$$
$${\rm and}~T_4 = \left(%
\begin{array}{ccccc}
  1 & 1 & 1 & 0 & 0 \\
  1 & 1 & 1 & 0 & 0 \\
  2 & 2 & 2 & 1 & 1 \\
  2 & 2 & 2 & 1 & 1 \\
  2 & 2 & 2 & 1 & 1 \\
\end{array}%
\right).$$

Now, by Lemma \ref{NewCondition}(4), we have
$${\displaystyle \sum_{i,j=1}^{5}t_{ij}^2=x(q+x)=60.}$$

However, the left-hand side of the last equality turns
out to be equal to 78 for $T_1$, 68 for $T_2$, 58 for $T_3$,
and 48 for $T_4$, a contradiction. This shows Theorem 1 for $x=6$.

The cases $x=4$ and $x=8$ can be drawn in exactly the same manner
(with more candidates for ${\cal T}(u)$ so we omit the details).

For $x=5$, the matrix
$$ \left(%
\begin{array}{ccccc}
  2 & 2 & 2 & 2 & 2 \\
  2 & 2 & 2 & 2 & 2 \\
  1 & 1 & 1 & 1 & 1 \\
  0 & 0 & 0 & 0 & 0 \\
  0 & 0 & 0 & 0 & 0 \\
\end{array}%
\right)$$
is the only pattern w.r.t. $u\in  {\overline {\cal L}}$ admissible by Lemma \ref{NewCondition},
while the matrix
$$ \left(%
\begin{array}{ccccc}
  4 & 4 & 2 & 2 & 2 \\
  4 & 4 & 2 & 2 & 2 \\
  2 & 2 & 0 & 0 & 0 \\
  2 & 2 & 0 & 0 & 0 \\
  2 & 2 & 0 & 0 & 0 \\
\end{array}%
\right)$$
is the only pattern w.r.t. $v\in  {\cal L}$.

Further, taking into account ${\cal T}(u)$, we see that there
exists a maximal clique $L$ of $G$ such that $|L\cap {\cal
L}|=10$, which contradicts ${\cal T}(v)$. This completes the proof
of Theorem 1. \wbull

\section{Uniqueness of Cameron -- Liebler line classes with $x=7$ in $PG(3,4)$}

Throughout the section we use the projective geometry notions. Let
us recall a construction of Cameron -- Liebler line class with
$x=7$ due to Govaerts and Penttila \cite{GovaertsPenttila}. We recall that
a {\it hyperoval} in a projective plane of order $q$ is a set of $q+2$ points,
 no 3 of which are collinear.

\smallskip

\noindent {\it Let $P$ be a point of $PG(3,4)$ and $\pi$ be a plane not containing $P$.
Let $O$ be a hyperoval in $\pi$ and $C$ be the set of lines incident to
the points of $O$ and $P$. Then $C$, all $2$-secants of $C$ and
all lines in $\pi$ external to $O$ form a Cameron -- Liebler line class
with parameter $x=7$. }
\smallskip

In this section we show the uniqueness of Cameron -- Liebler line class with $x=7$ in $PG(3,4)$
in the sense that any such line class arises in that way.
The proof essentially relies on the analysis of the patterns admissible by Lemma \ref{NewCondition}.

\begin{theo}
A Cameron -- Liebler line class in $PG(3,4)$ with $x=7$ is unique.
\end{theo}

The result is proven in several lemmas. In the remainder of this section,
let ${\cal L}$ be a Cameron -- Liebler line class with $x=7$ in $PG(3,4)$.

\begin{lemma} \label{7.0}
A line from ${\cal L}$ has one of the following patterns:
\begin{equation}\label{W.1}\tag{${\cal L}.1$}
\left(%
\begin{array}{ccccc}
  0 & 0 & 0 & 0 & 0 \\
  1 & 1 & 1 & 1 & 1 \\
  3 & 3 & 3 & 3 & 3 \\
  3 & 3 & 3 & 3 & 3 \\
  3 & 3 & 3 & 3 & 3 \\
\end{array}%
\right),
\end{equation}
\begin{equation}\label{W.2}\tag{${\cal L}.2$}
\left(%
\begin{array}{ccccc}
  4 & 4 & 2 & 3 & 2 \\
  4 & 4 & 2 & 3 & 2 \\
  3 & 3 & 1 & 2 & 1 \\
  2 & 2 & 0 & 1 & 0 \\
  2 & 2 & 0 & 1 & 0 \\
\end{array}%
\right),
\end{equation}
\begin{equation}\label{W.3}\tag{${\cal L}.3$}
\left(%
\begin{array}{ccccc}
  1 & 1 & 1 & 1 & 1 \\
  1 & 1 & 1 & 1 & 1 \\
  1 & 1 & 1 & 1 & 1 \\
  3 & 3 & 3 & 3 & 3 \\
  4 & 4 & 4 & 4 & 4 \\
\end{array}%
\right),
\end{equation}
while a line from ${\overline {\cal L}}$ has one of the following patterns:
\begin{equation}\label{B.1}\tag{${\overline {\cal L}}.1$}
\left(%
\begin{array}{ccccc}
  1 & 0 & 0 & 0 & 0 \\
  4 & 3 & 3 & 3 & 3 \\
  2& 1 & 1 & 1 & 1 \\
   2& 1 & 1 & 1 & 1 \\
   2& 1 & 1 & 1 & 1 \\
\end{array}%
\right),
\end{equation}
\begin{equation}\label{B.2}\tag{${\overline {\cal L}}.2$}
\left(%
\begin{array}{ccccc}
  1 & 0 & 0 & 0 & 0 \\
  1 & 0 & 0 & 0 & 0 \\
  3& 2 & 2 & 2 & 2 \\
  3& 2 & 2 & 2 & 2 \\
  3& 2 & 2 & 2 & 2 \\
\end{array}%
\right).
\end{equation}
\end{lemma}
\bproof The exhaustive enumeration of $5\times 5$-matrices admissible by Lemma \ref{NewCondition}.
\wbull

\begin{lemma}\label{BlockSet}
Let $L$ be a clique, such that $x<|C\cap {\cal L}|\leqslant x+q$. Then
$C\cap {\cal L}$ forms a blocking set in the plane, corresponding
to $L$. If there exist no Cameron--Liebler line classes with
parameter $x-1$, then the blocking set is nontrivial.
\end{lemma}
\bproof Lemma 4.1 (1) \cite{GovaertsPenttila}. \wbull

\begin{lemma} \label{7.1}
A line from ${\cal L}$ cannot have pattern $($\ref{W.3}$)$.
\end{lemma}
\bproof On the contrary, for a line that has pattern (\ref{W.3}),
we may relate any column of (\ref{W.3}) to a projective plane
with 11 points corresponding to lines from ${\cal L}$.
By Lemma \ref{BlockSet} and Theorem \ref{NonExistence},
we see that these points form a non-trivial blocking set.
However, it is easily seen from the last row of the pattern (\ref{W.3})
that the blocking set contains a line, a contradiction.\wbull

\begin{lemma} \label{7.2}
There exists a line that has pattern $($\ref{W.1}$)$.
\end{lemma}
\bproof It follows from the first row of patterns (\ref{B.1}) and (\ref{B.2})
that there exists a clique in $PG(3,4)$ containing a line $v$ from ${\cal L}$
and 20 lines from ${\overline {\cal L}}$.
Clearly, the line $v$ cannot have pattern (\ref{W.2}).
The lemma is proved.
\wbull

\begin{lemma}\label{7.3}
There exist six lines $l_{1},\ldots,l_{6}$ from ${\cal L}$ such that:

$(1)$ there is a point $P$ such that ${\rm Star}(P)\cap {\cal
L}=\{l_{1},\ldots,l_{6}\}$.

$(2)$ each of the lines $l_{1},\ldots,l_{6}$ has pattern $($\ref{W.1}$)$.

$(3)$ any plane contains either $0$ or $2$ lines from $l_{1},\ldots,l_{6}$.
\end{lemma}
\bproof Let a line $l_1$ have pattern (\ref{W.1}). The pattern
(\ref{W.1}) has the all-ones row. Therefore, up to the duality in
$PG(3,q)$, there exists a point $P$ such that ${\rm Star}(P)\cap {\cal
L}$ consists of six lines $l_{1},\ldots,l_{6}$, and any plane
containing $l_1$ has exactly one line among $l_{2},\ldots,l_{6}$.
This proves (1).

Suppose that the pattern w.r.t. the line $l_{2}$ is (\ref{W.2}).
Then one of the 3rd or 5th columns or of the last two rows of
(\ref{W.2}) corresponds to a clique in $PG(3,4)$ that contains the
lines $l_{1},\ldots,l_{6}$. This yields, without loss of
generality, that there exist a plane containing $l_{2}, l_{3},
l_{4}$, a plane containing $l_{2}, l_{5}, l_{6}$, and a plane
containing $l_{2}, l_{1}$. Let us consider the line $l_{3}$. The
pattern w.r.t. $l_3$ is (\ref{W.2}), since there is a pencil
containing three lines of ${\cal L}$ (namely, $l_{2}$, $l_{3},
l_{4}$). Therefore there exists one more plane, containing $l_3$
and a pair of lines chosen from $\{l_{1}, l_{5}, l_{6}\}$, a
contradiction. This gives (2).

Finally, Statement (3) follows from (1) and (2). The lemma is proved. \wbull
\medskip

Note that (\ref{W.1}) has the all-zero row. This means that a line
with pattern (\ref{W.1}) has a unique point incident to the only line of ${\cal L}$.
We call such a point {\it poor}.

\begin{lemma}\label{7.4} The following holds.

$(1)$ There exists a plane $\pi$ that contains all the poor points of lines $l_{1},\ldots,l_{6}$.

$(2)$ The poor points of lines $l_{1},\ldots,l_{6}$ form a hyperoval $O$ in $\pi$.

$(3)$ The $2$-secants of the set of poor points are the only lines from ${\overline {\cal L}}\cap {\rm line}(\pi)$.
\end{lemma}
\bproof Let $\pi$ be a plane containing three poor points of the
lines $l_1,l_2,l_3$, and a point $R$ of $l_6$. Suppose $R$ is not
poor. Since $R\ne P$ (recall that ${\rm Star}(P)\cap {\cal
L}=\{l_1,\ldots,l_6\}$), the point $R$ corresponds to the last
three rows of (\ref{W.1}), which is the pattern w.r.t. $l_6$.
Hence, $|{\rm Star}(R)\cap {\cal L}|=16$.

Now any line from ${\rm Star}(R)$ has pattern (\ref{W.1}), (\ref{W.2}), or (\ref{B.1})
(as only these patterns admit a clique in $PG(3,4)$ with 16 lines from ${\cal L}$). 
Further, it follows from those patterns that any pencil in a clique in $PG(3,4)$ with 16 lines
from ${\cal L}$ contains at least 3 lines from ${\cal L}$.
In other words, a pencil ${\rm pen}(R,\pi')$ contains at least 3 lines from ${\cal L}$
for any plane $\pi'\ni R$.

Let us show that there exists a line $l'\in {\rm pen}(R,\pi)$ that
is incident to at least two poor points of $l_1,l_2,l_3$. On the
contrary, suppose that each of the lines in ${\rm pen}(R,\pi)$ is
incident to at most one poor point of $l_1,l_2,l_3$. Then $R$ is
incident to at least three lines from ${\overline {\cal L}}$
(these are the lines, each of them is incident to a poor point of
$l_1,l_2,l_3$) and to at least three lines of ${\cal L}$ as $|{\rm
pen}(R,\pi)\cap {\cal L}|\geqslant 3$, a contradiction. 
Therefore, the line $l'$ exists. 

However, clearly, $l'$ has pattern (\ref{B.2}), which is impossible, and Statement (1) follows.

Statement (2) follows directly from Lemma \ref{7.3}(3).

Any 2-secant of the set $O$ of poor points of $l_1,\ldots,l_6$ has pattern (\ref{B.2}),
and the plane $\pi$ corresponds to one of the last four columns of (\ref{B.2}).
Hence, $\pi$ contains exactly 6 lines from ${\cal L}$ and $15={6 \choose 2}$ lines
from ${\overline {\cal L}}$, each of them is a 2-secant of $O$. The lemma is proved. \wbull

\begin{lemma}\label{7.5}
Any $2$-secant of lines $l_{1},\ldots,l_{6}$ external to poor points belongs to ${\cal L}$.
\end{lemma}
\bproof By Lemma \ref{Definitions}, a line of ${\cal L}$ meets exactly
50 lines of ${\cal L}$. By Lemma \ref{7.3}(3), a plane $\pi'\ni l_1$
contains exactly one more line of $l_{2},\ldots,l_{6}$, for example, $l_2$.
Now the line $l_1$ meets 10 lines of ${\overline {\cal L}}\cap {\rm line}(\pi')$:
the 2-secant of poor points of $l_1,l_2$, the six 1-secants of poor points of $l_1,l_2$,
and the three lines from ${\rm Star}(P)\cap {\rm line}(\pi')\setminus \{l_1,l_2\}$.
Therefore, the 9 remaining lines of $\pi'$ (i.e., 2-secants of $l_1$ and $l_2$ external
to their poor points) belong to ${\overline {\cal L}}$. The lemma is proved.\wbull
\medskip

Now, by Lemmas \ref{7.0}--\ref{7.5}, the Cameron -- Liebler line class ${\cal L}$
consists of the six lines $l_1,\ldots,l_6$ with common point $P$, the six lines
from $\pi$ external to the poor points of $l_1,\ldots,l_6$ that form a hyperoval $O$ in $\pi$,
and $9\times {6 \choose 2}=135$ lines that are 2-secants of pairs $l_i,l_j$, $i,j\in \{1,\ldots,6\}$.
Theorem 2 is proved.


\section{Cameron -- Liebler line classes in $PG(n,4)$}

The results of previous sections complete the classification of
Cameron -- Liebler line classes in $PG(3,4)$. The notion of
Cameron -- Liebler line class in $PG(3,q)$ can be naturally
generalized to that in $PG(n,q)$, see \cite{DrudgeThesis}. For a
subspace $X$ of $PG(n,q)$, we denote the set of all lines in $X$
by ${\rm line}(X)$.

\begin{lemma}\label{Definitions2} $(\cite[Theorem~3.2]{DrudgeThesis})$
Let ${\cal L}$ be a non-empty set of lines in $PG(n,q)$
with characteristic function $\chi_{\cal L}$. Then
the following are equivalent.

\quad $(1)$ There exists $x$ such that, for any flag $(P,X)$ with
$X$ an $i$-dimensional subspace of $PG(n,q)$, we have
$$|{\rm Star}(P)\cap {\cal L}| + \frac{\theta_{n-2}}{\theta_{i-1}\theta_{i-2}}|{\rm line}(X)\cap {\cal L}| = x + \frac{\theta_{n-2}}{\theta_{i-2}}|{\rm pen}(P,X)\cap {\cal L}|,$$
where $\theta_d:=q^d+\ldots+1$.

\quad $(2)$ There exists $x$ such that, for every line $l$,
$$|\{m\in {\cal L}\setminus \{l\}:~m~meets~l\}| = (q+1)x + (q^{n-1}+\ldots+q^2-1)\chi_{\cal L}(l).$$

In addition, if $n$ is odd:

\quad $(3)$ There exists $x$ such that $|{\cal S}\cap {\cal L}|=x$ for any line-spread ${\cal S}$.
\end{lemma}

A set ${\cal L}$ of lines of $PG(n,q)$ satisfying the above
conditions is called a {\it Cameron -- Liebler line class} in
$PG(n,q)$. In case of $n=3$ this definition coincides with that
from Section 2. The following are examples of Cameron -- Liebler
line classes:
\smallskip

-- the set of all lines in $PG(n,q)$,
\smallskip

-- the set of all lines contained in a hyperplane of $PG(n,q)$,
\smallskip

-- the set of all lines through a point,
\smallskip

-- for a non-incident point -- hyperplane pair $(P,H)$, the set of all lines through $P$ or in $H$.
\smallskip

Note that the complement to a Cameron -- Liebler line class with parameter $x$
is a Cameron -- Liebler line class with parameter $\frac{q^n+\ldots+1}{q+1}-x$.

Drudge \cite{DrudgeThesis} formulated the following generalized conjecture:
\smallskip

\noindent {\it The only Cameron -- Liebler line classes in $PG(n,q)$ are those
listed above and their complementary sets.}
\smallskip

In his Ph.D. Thesis \cite{DrudgeThesis}, Drudge proved this
conjecture for the case $q=3$. His approach essentially relies on
the classification of Cameron -- Liebler line classes in
$PG(3,3)$, see \cite{Drudge}, and, probably, can be applied to an
arbitrary $q$ once a classification of those in $PG(3,q)$ has been
obtained.

In this section, we prove the conjecture for $q=4$. For the
convenience of the reader, we give all the lemmas due to the
Drudge approach (for the proofs we refer to his Ph.D. Thesis). In
what follows, let ${\cal L}$ be a Cameron -- Liebler line class in
$PG(n,q)$, $n\geqslant 4$, $X$ be a three-dimensional subspace of
$PG(n,q)$.

\begin{lemma}\label{Reduction} $(\cite[Theorem~6.1]{DrudgeThesis})$
The line class ${\rm line}(X)\cap {\cal L}$ is a Cameron -- Liebler line class in $X$.
\end{lemma}

The intersection of ${\cal L}$ with $X$ can, in some cases, determine ${\cal L}$.

\begin{lemma}\label{StarIntersection} $(\cite[Lemma~6.1]{DrudgeThesis})$
If ${\rm line}(X)\cap {\cal L}$ consists of all the lines of $X$ on some point $P\in X$, i.e.,
${\rm line}(X)\cap {\cal L}={\rm Star}(P)\cap {\rm line}(X)$, then ${\cal L}={\rm Star}(P)$.
\end{lemma}

The lack of point -- plane duality in $PG(n,q)$, $n>3$, does not
allow to state immediately the similar lemma when
${\rm line}(X)\cap {\cal L}={\rm line}(\pi)$ for some plane $\pi$ of $X$.

\begin{lemma}\label{XtwoIntersection} $(\cite[Lemma~6.3]{DrudgeThesis})$
If ${\rm line}(X)\cap {\cal L}$ is a Cameron -- Liebler line class with parameter $x=2$ $($in $PG(3,q)$$)$
then ${\cal L}$ consists of all lines on a point plus all lines in a hyperplane,
for a non-icident point -- hyperplane pair in $PG(n,q)$.
\end{lemma}

A key step in the Drudge approach is a proof that a counterexample to the Cameron -- Liebler
conjecture on line classes in $PG(3,q)$ cannot occur as the intersection
of a three-dimensional subspace of $PG(n,q)$ with a Cameron -- Liebler line class of $PG(n,q)$.

\begin{lemma}\label{KeyLemma}
For $q=4$, the set ${\rm line}(X)\cap {\cal L}$ cannot be a
Cameron -- Liebler line class in $PG(3,4)$ with parameter $x=7$.
\end{lemma}
\bproof By taking complements if necessary, we may assume that
the set ${\rm line}(X)\cap {\cal L}$ is a Cameron -- Liebler line
class in $PG(3,4)$ with parameter $x=7$, which is unique by Theorem 2.
It follows from Lemma \ref{7.0} that there exists an incident
point--plane pair $(P,\pi)$ in $X$ such that $|{\rm line}(\pi)\cap {\cal L}|=11$,
$|{\rm pen}(P,\pi)\cap {\cal L}|=1$ and $|{\rm pen}(P,X)\cap {\cal L}|=1$
(we called such point {\it poor} in the previous section).


Let us show that ${\rm Star}(P)\cap {\cal L}$ consists of the only line,
which belongs to $\pi$. Clearly, $\pi$ and a line through
$P$ are contained in some three-dimensional                 
subspace $Y$ of $PG(n,4)$. By Lemma \ref{Reduction}, we see that
${\rm line}(Y)\cap {\cal L}$ is a Cameron -- Liebler line class in $Y$,
which intersects the plane $\pi\subset Y$ in 11 points.
However, the Cameron -- Liebler line class in $PG(3,4)$ with parameter $x=7$ is
the only line class that intersects a plane in 11 lines, see Lemma \ref{7.0}. 
Therefore ${\rm line}(Y)\cap {\cal L}$ is the Cameron -- Liebler
line class with parameter $x=7$, and, hence, ${\rm pen}(P,Y)\cap
{\cal L}={\rm pen}(P,X)\cap {\cal L}$. This yields that ${\rm
Star}(P)\cap {\cal L}={\rm pen}(P,X)\cap {\cal L}$, and every
line from ${\rm Star}(P)\setminus \pi$ does not belong to ${\cal L}$.

Further, since $n\geqslant 4$, there exists a three-dimensional subspace $Y$ on $P$
such that all lines of $Y$ on $P$ are not in ${\cal L}$ (note that
$\pi\not\subset Y$ in this case). By Lemma \ref{Reduction} and Theorem 2,
we see that a line class ${\rm line}(Y)\setminus {\cal L}$, which contains at least
all the lines of $Y$ on $P$, is a Cameron -- Liebler line class in $Y$ and its parameter $x\neq 7$.
By Theorem 1, this implies that
${\rm line}(Y)\setminus {\cal L}$ can only be one of the following line classes:
${\rm pen}(P,Y)$, 
${\rm line}(Y)$,  
${\rm pen}(P,Y)\cup {\rm line}(\pi')$ (where $P\notin \pi'$), 
or ${\rm line}(Y)\setminus {\rm line}(\pi')$ for some plane $\pi'\ne \pi$. 
Therefore $|{\rm line}(Y)\cap {\cal L}|$ equals $(17-1)\cdot 21$, $0$,
$(17-2)\cdot 21$, or $(17-1)\cdot 21$, respectively, whereas
$|{\rm pen}(P,Y)\cap {\cal L}|$ equals 0.

However, recall that $|{\rm Star}(P)\cap {\cal L}|=1$,
$|{\rm line}(X)\cap {\cal L}|=7\cdot (4^2+4+1)=147$, and
$|{\rm pen}(P,X)\cap {\cal L}|=1$.
Now applying Lemma \ref{Definitions2}~(1) to the subspaces $X$ and $Y$ gives
$$|{\rm Star}(P)\cap {\cal L}| + \frac{\theta_{n-2}}{\theta_{2}\theta_{1}}|{\rm line}(X)\cap {\cal L}| =
x + \frac{\theta_{n-2}}{\theta_{1}}|{\rm pen}(P,X)\cap {\cal L}|,$$
$$|{\rm Star}(P)\cap {\cal L}| + \frac{\theta_{n-2}}{\theta_{2}\theta_{1}}|{\rm line}(Y)\cap {\cal L}| =
x + \frac{\theta_{n-2}}{\theta_{1}}|{\rm pen}(P,Y)\cap {\cal L}|,$$ which is impossible. 
The lemma is proved.\wbull

%
%

\begin{theo}
For $q=4$, if a Cameron -- Liebler line class ${\cal L}$ is not of the type ${\rm Star}(P)$
for a point $P$ or ${\rm Star}(P)\cup {\rm line}(H)$ for a non-incident point -- plane pair $(P,H)$,
or the complement of a set of one of these two types then ${\cal L}$ is of the type ${\rm line}(H)$,
for a hyperplane $H$, or its complement.
In other words, a Cameron -- Liebler line class in $PG(n,4)$, $n\geqslant 4$,
is one of those listed above or their complementary sets.
\end{theo}
\bproof It follows from Theorem 1 and Lemmas
\ref{StarIntersection}, \ref{XtwoIntersection}, \ref{KeyLemma}
that a Cameron -- Liebler line class ${\cal L}$ satisfying the
hypothesis of the theorem intersects $X$ in one of the following two
ways:

\quad $(1)$ All lines of $X$, or no lines of $X$, or

\quad $(2)$ All lines in a plane of $X$, or all lines not in a plane of $X$.

Further, Lemma 6.5 and Theorem 6.3 from \cite{DrudgeThesis} show
that such a line class is of the type ${\rm line}(H)$, for a
hyperplane $H$, or its complement (in fact, these lemma and
theorem are proved under a slightly different assumption on
Cameron -- Liebler line classes in $PG(3,q)$, but it does not
matter in our case due to Lemma \ref{KeyLemma}). The theorem is
proved.\wbull
\medskip


We close this paper with several remarks.

Suppose we are given the set $G(v)\cap {\cal L}$ for some vertex $v$.
Then we are able to ``reconstruct'' the whole set ${\cal L}$ by using,
for instance, Lemma \ref{DefinitionsInGraphs}(\ref{Mu}). Actually,
the same idea was exploited in the study of completely regular codes
in the Johnson graphs \cite{GavGor}, and it led us to the study of those
in the Grassmann graphs.
\medskip

Recently Bamberg \cite{Bamberg} announced the non-existence of
Cameron -- Liebler line classes with parameter $x\in \{6,8\}$
in $PG(3,4)$ as a negative result of a computer search.
\medskip

We are grateful to Fr\'{e}d\'{e}ric Vanhove for the useful references.

\end{document}